\documentclass[runningheads]{svjour2}
\smartqed  % flush right qed marks, e.g. at end of proof
\journalname{Numerische Mathematik}
\usepackage{amssymb,amsmath}
\usepackage[dvips]{color}
\usepackage[dvips]{graphics}
\newtheorem{algorithm}{Algorithm}
\newtheorem{thm}{Theorem}
\newtheorem{lem}{Lemma}
\def\mat#1#2{\mathbb{R}^{#1\times #2}}
\newcommand{\bear}{\left(\begin{array}}
\newcommand{\enar}{\end{array}\right)}
\newcommand{\hot}{+ O(\macheps^2 )}
\newcommand{\beeq}{\begin{equation}}
\newcommand{\eneq}{\end{equation}}
\newcommand{\rank}{\mathrm{rank}}
\def\normtwo#1{\|#1\|_2}
\def\normF#1{\|#1\|_F}
\def\normal#1{#1^T\!~#1}
\def\ysqrt#1{\left( #1 \right)^{1/2}}
\def\co#1#2{#1\colon #2}
\newcommand{\mbf}{\mathbf}     %{\textbf{!!!}}
\newcommand{\tbf}{\textbf}
\def\eqref#1{{\normalfont(\ref{#1})}}
\def\numlist#1#2#3{#1=#2, \ldots , #3}
\newcommand{\macheps}{\varepsilon_M }
\def\matgen#1#2{(\mbf{#1}_1, \ldots , \mbf{#1}_{#2} )}

\begin{document}
\title{A note on the error analysis of classical {Gram--Schmidt}}
\author{Alicja Smoktunowicz \and Jesse L. Barlow
%\thanks{Jesse Barlow's research was supported by the National Science Foundation under  grant no.
%CCF-0429481.}
\and Julien Langou}
 \institute{
Alicja Smoktunowicz \at
 Faculty of Mathematics and
Information Science, Warsaw University of Technology, Pl.
Politechniki 1, Warsaw, 00-661 Poland, \email{smok@mini.pw.edu.pl}
\and
Jesse L. Barlow \at
 Department of Computer Science and Engineering, The
Pennsylvania State University, University Park, PA 16802-6822,
USA, \email{barlow@cse.psu.edu} \and
Julien Langou \at Department of Computer
Science, The University of Tennessee, 1122 Volunteer Blvd.,
Knoxville, TN 37996-3450, USA, \email{langou@cs.utk.edu}
}

\date{   }
%\date{Received: date / Revised version: date}
\maketitle

\begin{abstract}
An error analysis result is given  for classical Gram--Schmidt factorization of a full rank matrix $A$ into $A=QR$
where $Q$ is left orthogonal (has orthonormal columns) and $R$ is upper triangular. The work presented here shows that
the computed $R$  satisfies $\normal{R}=\normal{A}+E$ where $E$ is an appropriately small backward error, but only if the diagonals of $R$ are computed in
a manner similar to Cholesky factorization of the normal equations matrix.

A similar result is stated  in [Giraud at al, Numer. Math. 101(1):87--100,2005]. However, for that result to hold, the diagonals
of $R$ must be computed in the manner recommended in this work.
\end{abstract}

\vspace{3ex}

The classical Gram--Schmidt (CGS) orthogonal factorization is analyzed in a
recent work of Giraud et al. \cite{Gir05} and in a number of other sources
\cite{Bjo67a,Kie74,Wol01,BSE05,dgks76,Hoff89}, \cite[\S 6.9]{Parl98}, \cite[\S
2.4.5]{Bjo96}.

For a matrix $A \in \mat{m}{n}$ ($m \geq n$) with $\rank(A)=n$, in
exact arithmetic, the algorithm produces a factorization \beeq A =
Q R \label{eq:QRfact} \eneq where $Q$ is {\em left orthogonal}
(i.e. $\normal{Q} = I_n$), and $R \in \mat{n}{n}$ is upper
triangular and nonsingular.  In describing the algorithms, we use
the notational conventions, \begin{eqnarray*}A &=&\matgen{a}{n},
\quad Q = \matgen{q}{n},\\
R&=& (r_{jk}). \end{eqnarray*}

The algorithm forms $Q$ and $R$ from $A$ column by column as described in the
following  pseudo-code. We label this algorithm CGS--S, for classical
Gram--Schmidt ``standard.''

%\newpage
\noindent
\begin{algorithm}[Classical Gram--Schmidt Orthogonal Factorization (Standard) (CGS--S)] $\:$
\label{alg:CGS}
\begin{tabbing}
$r_{11} = \normtwo{\mbf{a}_1};  \mbf{q}_1 = \mbf{a}_1/r_{11};$\\
$R_1=(r_{11}); Q_1 = (\mbf{q}_1); $ \\
\tbf{for} \= $k=\co{2}{n}$ \\
\>$\mbf{s}_k =Q_{k-1}^T \mbf{a}_k;$
\\
\>$\mbf{v}_k = \mbf{a}_k-Q_{k-1} \mbf{s}_k ;$
\\
\>${r}_{kk} = \normtwo{\mbf{v}_k};$
\\
\>$\mbf{q}_k= \mbf{v}_k/r_{kk};$
\\
\> $R_k =\bordermatrix{ & k-1 & 1 \cr k-1 & R_{k-1} & \mbf{s}_k
\cr 1 & 0 & r_{kk} } ;\,\,  Q_k = \bordermatrix{ & k-1 & 1 \cr &
Q_{k-1} &
\mbf{q}_k };$ \\
\tbf{end}; \\
$Q=Q_n; \,\, R=R_n$;\\
\end{tabbing}
\end{algorithm}

As is well known \cite[p.63,\S 2.4.5]{Bjo96}, in floating point arithmetic, $Q$
is far from left orthogonal. The authors of \cite{Gir05} prove a number of
results about classical Gram--Schmidt. This note shows that for one of their
results (Lemma 1 in \cite{Gir05}), the diagonal elements $r_{kk}$ should be
computed differently from  Algorithm \ref{alg:CGS}, substituting a
Cholesky-like formula for $r_{kk}$ rather than setting  $r_{kk}
=\normtwo{\mbf{v}_k}$.  That change produces the Algorithm \ref{alg:SmokBar}.
Since it uses a pythagorean identity to compute the diagonals of $R$, we call
it CGS-P for ``classical Gram--Schmidt pythagorean.''.

%\newpage
\noindent
\begin{algorithm}[Cholesky--like Classical Gram--Schmidt Orthogonal Factorization (CGS--P)] $\:$
\label{alg:SmokBar}
\begin{tabbing}
$r_{11} = \normtwo{\mbf{a}_1};  \mbf{q}_1 = \mbf{a}_1/r_{11};$\\
$R_1=(r_{11}); Q_1 = (\mbf{q}_1); $ \\
\tbf{for} \= $k=\co{2}{n}$ \\
\>$\mbf{s}_k =Q_{k-1}^T \mbf{a}_k;$
\\
\>$\mbf{v}_k = \mbf{a}_k-Q_{k-1} \mbf{s}_k ;$
\\
\> $\psi_k = \normtwo{\mbf{a}_k}; \phi_k = \normtwo{\mbf{s}_k}$; \\
\>${r}_{kk} = \ysqrt{\psi_k-\phi_k} \ysqrt{\psi_k+\phi_k}$;
\\
\>$\mbf{q}_k= \mbf{v}_k/r_{kk};$
\\
\> $R_k =\bordermatrix{ & k-1 & 1 \cr k-1 & R_{k-1} & \mbf{s}_k
\cr 1 & 0 & r_{kk} } ; \,\, Q_k = \bordermatrix{ & k-1 & 1 \cr &
Q_{k-1} &
\mbf{q}_k };$ \\
\tbf{end}; \\
$Q=Q_n; R=R_n$;\\
\end{tabbing}
\end{algorithm}

%This computation of $r_{kk}$ does not affect the results in
%\cite{Gir05} on re-orthogonalization of $\mbf{q}_k$ against
%$Q_{k-1}$.
%Our numerical tests indicate that  the main and quite useful result  in \cite[Theorem
%2]{Gir05}  still holds if $r_{kk}$ is computed as in Algorithm \ref{alg:CGS}
% {\em after the re-orthogonalization}.

We assume that we are using a floating point arithmetic that satisfies the
IEEE floating point standard. In IEEE   arithmetic
\[
f\ell(x+y) = (x+y) (1+ \delta), \quad |\delta| \leq \macheps
\]
for results in the normalized range \cite[p.32]{Ove01}.

Letting $\macheps$ be the machine unit, we follow Golub and Van Loan \cite[\S 2.4.6]{GoVa96} and use the linear approximation
\[
(1+\macheps)^{p(n)} = 1+p(n) \macheps \hot
\]
for a modest function $p(n)$ thereby assuming that the $O(\macheps^2)$  makes no significant contribution.

For the sake of self containment, we give Lemma 1 from \cite{Gir05}.

\begin{lemma} \cite{Gir05} In floating point arithmetic with machine unit $\macheps$, the computed upper triangular factor from Algorithm \ref{alg:CGS}
satisfies
\[
\normal{R} = \normal{A}+E, \quad \normtwo{E} \leq c(m,n) \normtwo{A}^2 \macheps
\]
where $c(m,n)=O(mn^2)$.
\label{lem:Giretal}
\end{lemma}

As stated, this lemma is not correct for Algorithm \ref{alg:CGS}, but a slightly different version of this result holds for Algorithm \ref{alg:SmokBar}.

 We define the four functions
\begin{eqnarray}
c_1(m,k)&=& \left\{ \begin{array}{ll} 1 & k=1 \\ 2\sqrt{2}mk+2\sqrt{k}& \numlist{k}{2}{n}, \end{array} \right.\nonumber \\
c_2(m,k)&=& \left\{ \begin{array}{ll} m+2 & k=1 \\ 3.5mk^2- 1.5mk+16k& \numlist{k}{2}{n}, \end{array} \right. \label{eq:c2mk}\\
c_3(m,k)&=& 0.5 c_2(m,k), \quad c_4(m,k)= c_2(m,k)+2c_1(m,k), \nonumber
\end{eqnarray}
we let $A_k$ be the first $k$ columns of $A$, and let
\[
\kappa_2(R_k) = \normtwo{R_k} \normtwo{R_k^{-1}}.
\]
The new version of Lemma \ref{lem:Giretal} is
Theorem \ref{thm:stepk}.

\begin{thm}
\label{thm:stepk} Assume that in floating point arithmetic with
machine unit $\macheps$, for the $R$ resulting from Algorithm
\ref{alg:SmokBar} for each $k$, we have \beeq c_4(m,k) \macheps \kappa_2(R_k)^2
 < 1. \label{eq:mainassumption} \eneq
Let $A_k \in \mat{m}{k}$ consist of the first $k$ columns of $A$.
 Then, for
$\numlist{k}{1}{n}$, to within terms of $O(\macheps^2)$, the computed matrices $R_k$ and $Q_k$ satisfy
\begin{eqnarray}
 Q_k R_k-A_k &=& \Delta A_k,  \quad \normtwo{\Delta A_k } \leq
c_1(m,k) \normtwo{A_k} \macheps , \label{eq:backerr}
\\
\normal{R_k}-\normal{A_k} &=&  E_k, \quad \normtwo{E_k} \leq
c_2(m,k) {\normtwo{A_k}}^2 \macheps, \label{eq:normaleqbnd}
\\
\normtwo{R_k}&=& \normtwo{A_k} (1+\mu_k), \quad |\mu_k| \leq
c_3(m,k) \macheps, \label{eq:normRkbnd}
\\
\normtwo{I - \normal{Q_k}} &\leq& c_4(m,k) \kappa_2(R_k)^2
 \macheps, \label{eq:Qkbnd}
\\
\normtwo{Q_k} &\leq& \sqrt{2}. \label{eq:normQkbnd}
\end{eqnarray}
\end{thm}

The proof of Theorem \ref{thm:stepk} is given in the appendix.

The restriction (\ref{eq:mainassumption}) assures that $R$ is nonsingular, and that (\ref{eq:Qkbnd}) and (\ref{eq:normQkbnd}) hold.  A weaker assumption that assures that $R$ is nonsingular and that $\normtwo{Q_k}$ is bounded would
yield bounds similar to (\ref{eq:backerr}), (\ref{eq:normaleqbnd}), and (\ref{eq:normRkbnd}).

\begin{remark}
\label{rem:condition}
The condition (\ref{eq:mainassumption}) and the bound (\ref{eq:Qkbnd}) are  stated in terms
of $\kappa_2(R_k)$. We now show how it may be stated in terms of
\[
\kappa_2(A_k) = \normtwo{A_k} \normtwo{A_k^{\dagger}}
\]
where $A_k^{\dagger}$ is the Moore-Penrose pseudoinverse of $A_k$. In exact arithmetic, $\kappa_2(A_k)$ and
$\kappa_2(R_k)$ are the same quantity, and equation (\ref{eq:normRkbnd}) states that $\normtwo{R_k}$ and $\normtwo{A_k}$ are
nearly interchangable in floating point arithmetic. To relate $\normtwo{R_k^{-1}}$ and $\normtwo{A_k^{\dagger}}$, we use
 eigenvalue inequalities.

From the fact that
\beeq
\normtwo{R_k^{-1}}^{-1}= \sqrt{\lambda_k(\normal{R_k})}, \quad \normtwo{A_k^{\dagger}}^{-1}= \sqrt{\lambda_k(\normal{A_k})}
\label{eq:lambdakdef}
\eneq
where $\lambda_k(\cdot)$ denotes $kth$ largest (and therefore smallest) eigenvalue, we can obtain an upper bound for $\normtwo{A_k^{\dagger}}$ using
Weyl's monotonicity theorem \cite[Theorem 10.3.1]{Parl98}. Applying that theorem to (\ref{eq:normaleqbnd}), we have
\begin{eqnarray*}
\lambda_k(\normal{R_k}) &\geq& \lambda_k(\normal{A_k}) -\normtwo{E_k} \\
&\geq&  \lambda_k(\normal{A_k}) -\macheps c_2(m,k) \normtwo{A_k}^2 \hot \\
&=& \lambda_k(\normal{A_k}) -\macheps c_2(m,k) \normtwo{R_k}^2 \hot \\
&\geq& \lambda_k(\normal{A_k})(1-\zeta_k)
\end{eqnarray*}
where
\beeq
\zeta_k =\macheps c_2(m,k) \kappa_2(R_k)^2 \hot.
\label{eq:zetakdef}
\eneq
Using (\ref{eq:lambdakdef}), we have
\[
\normtwo{R_k^{\dagger}} \leq \normtwo{A_k^{-1}}(1-\zeta_k)^{-1/2}.
\]
From (\ref{eq:normRkbnd}), we may conclude that
\[
\kappa_2(R_k) \leq \kappa_2(A_k) (1+\mu_k)(1-\zeta_k)^{-1/2}.
\]
Thus a slight variation of the condition (\ref{eq:mainassumption}) may be stated in terms of $\kappa_2(A_k)$.
 Since it fits more naturally into the proof of Theorem \ref{thm:stepk} and it is more
easily computed than $\kappa_2(A_k)$, we use $\kappa_2(R_k)$.
\end{remark}

The conclusion of Theorem \ref{thm:stepk} does not hold for Algorithm \ref{alg:CGS}, as shown by the following example.
We were able to construct several similar examples.  Both examples were done in MATLAB version 7 on a Dell Precision 370 workstation
running Linux.

\begin{example}
\label{ex:counterexample}
We produced a $6 \times 5$ matrix with the following MATLAB code.

\begin{tabbing}
B=hilb(6); \\
$A1=ones(6,3)+B(\colon,\co{1}{3})*1e-2$;\\
B=pascal(6);\\
$A2=B(\colon,\co{1}{2})$;\\
A=[A1 A2];\\
\end{tabbing}

The command hilb(6) produces the $6 \times 6$ Hilbert matrix, the
command ones(6,3) produces a $6 \times 3$ matrix of ones, and the
command pascal(6) produces a $6 \times 6$ matrix from Pascal's
triangle. The condition number of $R$ from  Algorithm 2,$\, \kappa_2(R) =
\normtwo{R}\normtwo{R^{-1}}$, computed by the MATLAB command \tbf{cond}, is $3.9874\cdot 10^6$, thus given that
$\macheps \approx 2.2206 \cdot 10^{-16}$ in IEEE double precision,
$R$ is neither well-conditioned nor near singular.

We computed the Q--R factorization using Algorithm \ref{alg:CGS} (CGS--S)
and then we computed the same factorization using Algorithm \ref{alg:SmokBar} (CGS--P). The
resulting $Q$ and $R$ satisfy the results in Table \ref{tab:error}.

\begin{table}
\label{tab:error}

\begin{center}
\begin{tabular}{|c|c|c|c|}\hline
Algorithm & $\normtwo{\normal{A}-\normal{R}}/\normtwo{A}^2$ &
$\normtwo{I-\normal{Q}}$ \\ \hline CGS--S (Algorithm \ref{alg:CGS}) &  4.5460e-9
& 3.9874e-6 \\ \hline
CGS--P (Algorithm \ref{alg:SmokBar} &  3.3760e-17 & 5.2234e-5 \\
\hline
\end{tabular}
\end{center}
\caption{Orthogonality and Normal Equations Error from CGS Algorithms for Example \ref{ex:counterexample}}
\end{table}
\end{example}

The bound on $\normtwo{\normal{A}-\normal{R}}$ in
(\ref{eq:normaleqbnd}) appears to be satisfied if  $r_{kk}$ is
computed as in Algorithm \ref{alg:SmokBar}, but it is not if
$r_{kk}$ is computed as in  Algorithm \ref{alg:CGS}.

A larger, more complex, but better conditioned example is given next. 

\begin{example}
\label{ex:gluedmatrix}
A large class of examples where CGS-S obtains a large value of $\normtwo{\normal{A}-\normal{R}}/(\normtwo{A}^2)$, but
CGS-P arises from glued matrices. A general MATLAB code for these glued matrices is given by

\begin{verbatim}
function [A]=create_gluedmatrix (condA_glob,condA,m,nglued,nbglued)
	n = nglued*nbglued;
	A = orth(rand(m,n));
	A = A*diag([10.^(0:condA_glob/(n-1):condA_glob)])*orth(randn(n,n));
	ibeg = 1;
	iend = nglued;
	for i=1:nbglued,
		A(:,ibeg:iend) = A(:,ibeg:iend)*diag([10.^(0:condA/(nglued-1):condA)])...
                *orth(randn(nglued,nglued));
		ibeg = ibeg+nglued;
		iend = iend+nglued;
	end
\end{verbatim}
Here $m$ represents the number of rows of $A$, $nglued$ is the number of columns in a block, $nbglued$ is the number of
blocks that are glued together, and $n= nglued \times nbglued$ is the number of columns in the matrix. The parameter
$condA$ is the condition number of a block, and $condA\_glob$ is a parameter to couple the blocks together. The MATLAB command \verb+orth(X)+ 
produces an orthonormal basis for the range of $X$, thus the command \verb+orth(randn(m,n))+ produces a random orthogonal matrix.

For this example, we used the parameters
\[
condA\_glob=1; condA=2; m=200; nglued=5; nbglued=40;
\]
for which we obtained a $200 \times 200$ matrix with condition number $506.92$ (the condition number of the orthogonal factor $R$ is about the same).
 We also used the command \verb+randn('state',0)+ to reset the random number generator to its initial state. Table \ref{tab:glued} summarizes the results
from applying CGS--S and CGS--P to this matrix.

\begin{table}
\label{tab:glued}

\begin{center}
\begin{tabular}{|c|c|c|c|}\hline
Algorithm & $\normtwo{\normal{A}-\normal{R}}/\normtwo{A}^2$ &
$\normtwo{I-\normal{Q}}$ \\ \hline CGS--S (Algorithm \ref{alg:CGS}) &  3.8744e-6
& 9.3676e-4 \\ \hline
CGS--P (Algorithm \ref{alg:SmokBar}) &  2.8729e-16 & 1.8972e-12 \\
\hline
\end{tabular}
\end{center}
\caption{Orthogonality and Normal Equations Error from CGS Algorithms for Example \ref{ex:gluedmatrix}}
\end{table}

For this  example,  the loss of orthogonality of CGS--S is far in excess of $O(\epsilon \kappa_2(R)^2)$, whereas the loss of orthogonality for
CGS--P is well within that bound. The error $\normtwo{\normal{A}-\normal{R}}$ is far larger for CGS--S than it is for CGS--P and is much
greater than $O(\macheps \normtwo{A}^2)$. 
\end{example}

%For both
%methods of computing $r_{kk}$, the orthogonality error in $Q$ is
%about what would be expected from \cite[Theorem 1]{Gir05} and Theorem \ref{thm:stepk} above.

\section*{Conclusion} The upper triangular factor $R$ from classical Gram--Schmidt has been
shown to satisfy the bound (\ref{eq:normaleqbnd}) provided that the diagonal elements of $R$ are computed as they are in the
Cholesky factorization of the normal equations matrix. If these diagonal elements are computed as in standard versions of classical Gram--Schmidt,
no bounds  such as (\ref{eq:normaleqbnd}) or (\ref{eq:Qkbnd}) may be guaranteed.

%\bibliographystyle{spmpsci}
%\bibliography{stat,cuppen,proposals}

\section*{Appendix. Proof of Theorem \ref{thm:stepk}}

To set up the proof of Theorem \ref{thm:stepk}, we
require a lemma.

\begin{lem}
\label{lem:rkk}
Let $Q \in \mat{m}{n}$ and $R\in \mat{n}{n}$ be the results of Algorithm \ref{alg:SmokBar} in floating point arithmetic with machine unit
$\macheps$ and that $R$ satisfies \eqref{eq:mainassumption}. Then
\begin{equation}
 r_{11} = \normtwo{\mbf{a}_1} (1+ \delta_1), \quad |\delta_1|
\leq (0.5 m + 1) \macheps \hot \label{eq:point8}
\end{equation}
and for $k=2, \ldots, n$
\begin{eqnarray}
r_{kk}&=& \ysqrt{\normtwo{\mbf{a}_k}^2 (1+\delta_k) -
\normtwo{\mbf{s}_k}^2 (1+\Delta_k)},
\label{eq:point9}
\\
|\delta_k|, |\Delta_k| &\leq& (m +8) \macheps  \hot, \nonumber
\\
\normtwo{\mbf{s}_k} &\leq& \normtwo{\mbf{a}_k} (1+\zeta), \quad
|\zeta | \leq (m + 2) \macheps \hot. \label{eq:point10}
\end{eqnarray}
\end{lem}

\begin{proof}
Equation (\ref{eq:point8}) is just the error in the computation of
$\normtwo{\mbf{a}_1}$. In the computation of $r_{kk},
\numlist{k}{2}{n}$, note that
\begin{eqnarray}
\psi_k &=& f\ell(\normtwo{\mbf{a}_k}) = \normtwo{\mbf{a}_k} (1 + \epsilon_1^{(k)}), \label{eq:psiinfloat}
\\
\phi_k &=& f\ell(\normtwo{\mbf{s}_k}) = \normtwo{\mbf{s}_k} (1 + \epsilon_2^{(k)}), \label{eq:phiinfloat}
\\
|\epsilon_i^{(k)}| &\leq& (0.5 m +1)\macheps \hot, \quad i=1,2. \nonumber
\end{eqnarray}

  Using (\ref{eq:mainassumption}), we conclude that $R$ is nonsingular, thus $r_{kk} > 0$ for all $k$. Thus in Algorithm \ref{alg:SmokBar}, $r_{kk} >0$ only if $\psi_k > \phi_k$.

To get (\ref{eq:point9}), note that
\[
r_{kk} = \sqrt{\psi_k - \phi_k}\sqrt{\psi_k + \phi_k} (1+\epsilon_3^{(k)}), \quad |\epsilon_3^{(k)}| \leq 3\macheps \hot.
\]
Thus using (\ref{eq:psiinfloat}) and (\ref{eq:phiinfloat}), we have
\begin{eqnarray*}
r_{kk} &=& \sqrt{\normtwo{\mbf{a}_k}^2 (1+\epsilon_1^{(k)})^2-\normtwo{\mbf{s}_k}^2 (1+\epsilon_2^{(k)})^2}(1+\epsilon_3^{(k)})\\
&=&\ysqrt{\normtwo{\mbf{a}_k}^2 (1+\delta_k) -
\normtwo{\mbf{s}_k}^2 (1+\Delta_k)}
\end{eqnarray*}
where
\begin{eqnarray*}
\delta_k &=& (1+\epsilon_1^{(k)})^2 (1+\epsilon_3^{(k)})^2 -1, \\
\Delta_k &=& (1+\epsilon_2^{(k)})^2 (1+\epsilon_3^{(k)})^2 -1.
\end{eqnarray*}
That yields
\[
|\delta_k | , |\Delta_k | \leq (m+8)\macheps \hot.
\]
Therefore $r_{kk}$ satisfies (\ref{eq:point9}).

Since $\psi_k > \phi_k$ as outlined above, from (\ref{eq:psiinfloat})--(\ref{eq:phiinfloat}), we have
\[
\psi_k = \normtwo{\mbf{a}_k}(1 + \epsilon_1^{(k)}) > \phi_k = \normtwo{\mbf{s}_k} (1 + \epsilon_2^{(k)})
\]
thus
\begin{eqnarray*}
\normtwo{\mbf{s}_k } &<& \normtwo{\mbf{a}_k}(1 + \epsilon_1^{(k)})(1 + \epsilon_2^{(k)})^{-1}\\
& \leq& \normtwo{\mbf{a}_k} (1+\zeta)
\end{eqnarray*}
where $\zeta$ satisfies (\ref{eq:point10}).
\end{proof}

As a consequence of  the singular value version of the Cauchy interlace theorem \cite[p.449-450, Corollary 8.6.3]{GoVa96},
we have that $\normtwo{R_k} \leq \normtwo{R}$ and $\normtwo{R_k^{-1}} \leq \normtwo{R^{-1}}$. We will use these facts freely
in the proof of Theorem \ref{thm:stepk}.

We can now prove Theorem \ref{thm:stepk}.

\begin{proof}[of Theorem \ref{thm:stepk}]
The results (\ref{eq:backerr})--(\ref{eq:normaleqbnd}) are proven by induction on $k$.
First, consider $k=1$. From Lemma \ref{lem:rkk}, we have (\ref{eq:point8}), so
\[
r_{11} = \normtwo{\mbf{a}_1} (1+ \delta_1), \quad |\delta_1| \leq (0.5m+1) \macheps \hot
\]
which implies that
\begin{eqnarray*}
\normal{R_1} &=& r_{11}^2 = \normtwo{\mbf{a}_1}^2 (1+ \delta_1)^2 \\
&=& \normal{A_1}(1+\delta_1)^2 = \normal{A_1} +E_1
\end{eqnarray*}
where
\[
E_1 = 2\delta_1 \normal{A_1} + \delta_1^2 \normal{A_1}.
\]
Thus
\[
\normtwo{E_1} = |E_1| \leq (m + 2)  \normtwo{\mbf{a}_1}^2 \macheps
\hot = (m+2) \normtwo{A_1}^2 \macheps \hot.
\]
Also, we can conclude from standard error bounds that
\[
\mbf{q}_1 = (I+G_1)\mbf{a}_1/r_{11}, \quad \normtwo{G_1}\leq \macheps.
\]
Therefore
\[
A_1 -Q_1 R_1 = \mbf{a}_1 - \mbf{q}_1 r_{11} = -G_1 \mbf{a}_1
\]
so that
\beeq
\normtwo{A_1 -Q_1 R_1} = \normtwo{\mbf{a}_1 - \mbf{q}_1 r_{11}} \leq \normtwo{ G_1}\normtwo{ \mbf{a}_1} \leq \macheps \normtwo{ \mbf{a}_1}.
\label{eq:DeltaA1}
\eneq

Assume that  (\ref{eq:backerr})--(\ref{eq:normQkbnd}) hold for $k-1$, and prove them for $k$.  We first prove (\ref{eq:backerr})--(\ref{eq:normaleqbnd}), and then show that (\ref{eq:normRkbnd})--(\ref{eq:normQkbnd}) follow.

First, we start with error bounds of the computation of the vectors $\mbf{s}_k$,$\mbf{v}_k$, and $\mbf{q}_k$ to prove (\ref{eq:backerr}).
Note that
\beeq
\mbf{s}_k = f\ell(Q_{k-1}^T \mbf{a}_k) = Q_{k-1}^T \mbf{a}_k -\delta \mbf{s}_k
\label{eq:skbnd1}
\eneq
where
\begin{eqnarray}
\normtwo{\delta \mbf{s}_k} &\leq& m\sqrt{k-1}  \normtwo{Q_{k-1}} \normtwo{\mbf{a}_k}\macheps \hot \nonumber \\
&\leq& \sqrt{2(k-1)}m  \normtwo{\mbf{a}_k} \macheps \hot.
\label{eq:skbnd2}
\end{eqnarray}
Also, we have
\beeq
\mbf{v}_k = f\ell(\mbf{a}_k - Q_{k-1}\mbf{s}_k)  = \mbf{a}_k - Q_{k-1}\mbf{s}_k-\delta \mbf{v}_k
\label{eq:vkbnd1}
\eneq
where
\[
\normtwo{\delta\mbf{v}_k} \leq  \normtwo{\mbf{a}_k} \macheps + \sqrt{k-1}m \normtwo{Q_{k-1}} \normtwo{\mbf{s}_k}\macheps \hot.
\]
>From (\ref{eq:point10}), the bound on $\normtwo{\mbf{s}_k}$ in
(\ref{eq:point10}), and the induction hypothesis on $Q_{k-1}$, we
have \beeq \normtwo{\delta \mbf{v}_k} \leq
(\sqrt{2(k-1)}m+1)\normtwo{\mbf{a}_k}\macheps \hot.
\label{eq:vkbnd2} \eneq Again using the bound on
$\normtwo{\mbf{s}_k}$ in (\ref{eq:point10}), we note that
\begin{eqnarray*}
{\normtwo{\mbf{v}_k+\delta \mbf{v}_k}}^2 &=& \normtwo{\mbf{a}_k}^2 -2\mbf{a}_k^T Q_{k-1}\mbf{s}_k + \normtwo{Q_{k-1}\mbf{s}_k}^2 \\
&=&\normtwo{\mbf{a}_k}^2 -2\normtwo{\mbf{s}_k}^2 + \normtwo{Q_{k-1}\mbf{s}_k}^2-2(\delta\mbf{s}_k)^T \mbf{s}_k \\
&\leq&\normtwo{\mbf{a}_k}^2 -2\normtwo{\mbf{s}_k}^2 + \normtwo{Q_{k-1}}^2 \normtwo{\mbf{s}_k}^2-2(\delta\mbf{s}_k)^T \mbf{s}_k \\
&\leq& \normtwo{\mbf{a}_k}^2 -2\normtwo{\mbf{s}_k}^2 +2\normtwo{\mbf{s}_k}^2-2(\delta\mbf{s}_k)^T \mbf{s}_k \\
&=& \normtwo{\mbf{a}_k}^2 -2(\delta \mbf{s}_k)^T \mbf{s}_k \\
&\leq& \normtwo{\mbf{a}_k}^2 +2\normtwo{\delta \mbf{s}_k}\normtwo{\mbf{s}_k} \\
&=& \normtwo{\mbf{a}_k}^2 + 2 \normtwo{\delta \mbf{s}_k}\normtwo{\mbf{a}_k} \hot \\
&\leq& \normtwo{\mbf{a}_k}^2(1+ \sqrt{2(k-1)}m \macheps)^2 \hot.
\end{eqnarray*}
Thus
\[
\normtwo{\mbf{v}_k} \leq \normtwo{\mbf{a}_k}(1+(3\sqrt{2(k-1)}m)\macheps) \hot =\normtwo{\mbf{a}_k} +O(\macheps).
\]
We note that
\[
\mbf{q}_k = (I+G_k)\mbf{v}_k/r_{kk}, \quad \normtwo{G_k} \leq \macheps.
\]
If we let
\[
\Delta A_k = Q_k R_k - A_k
\]
then
\[
\Delta A_k = \bear{cc} \Delta A_{k-1} & \delta \mbf{a}_k \enar
\]
where
\begin{eqnarray*}
\delta \mbf{a}_k &=& (I+G_k) \mbf{v}_k + Q_{k-1} \mbf{s}_k -\mbf{a}_k, \\
&=& G_k \mbf{v}_k - \delta \mbf{v}_k.
\end{eqnarray*}
That yields
\[
\normtwo{\delta \mbf{a}_k} \leq \normtwo{G_k} \normtwo{\mbf{v}_k } + \normtwo{\delta \mbf{v}_k} \leq (2\sqrt{2(k-1)}m+2)\macheps \normtwo{\mbf{a}_k} \hot.
\]

To bound $\normtwo{\Delta A_k}$, we give a recurrence for bounding
$\normF{\Delta A_k}$ in terms of $\normF{A_k}$, then use the bound
$\normF{A_k} \leq \sqrt{k} \normtwo{A_k}$. We show that
\[
\normF{\Delta A_k} \leq \hat{c}_1 (m,k) \normF{A_k} \macheps \hot.
\]
For $k=1$,
\[
\normF{\Delta A_1} = \normtwo{\mbf{a}_1} = \macheps \normtwo{\mbf{a}_1} = \macheps \normF{A_1}.
\]
Using properties of the Frobenius norm,
\begin{eqnarray}
 \normF{\Delta A_k}^2 &\leq&  \normF{\Delta A_{k-1}}^2 +\normtwo{\delta \mbf{a}_k}^2 \nonumber \\
&\leq& [\hat{c}_1^2(m,k-1) \normF{A_{k-1}}^2 + (2\sqrt{2(k-1)}m+2)^2\normtwo{\mbf{a}_k}^2 ]\macheps^2  + O(\macheps^3) \nonumber \\
 &\leq&  \max\{ \hat{c}_1^2(m,k-1),(2\sqrt{2(k-1)}m+2)^2  \} (\normF{A_{k-1}}^2+\normtwo{\mbf{a}_k}^2)\macheps^2 + O(\macheps^3) \nonumber\\
&=& \hat{c}_1^2(m,k) \normF{A_k}^2 \macheps^2 + O(\macheps^3).
\label{eq:DeltaAkbnd}
\end{eqnarray}
A quick induction argument yields
\[
\hat{c}_1(m,k) = 2 \sqrt{2(k-1)}m+2 \leq 2\sqrt{2k}m+2.
\]
Thus
\[
\normtwo{\Delta A_k} \leq \normF{\Delta A_k} \leq \hat{c}_1 (m,k) \macheps \normF{A_k} \hot \leq \sqrt{k} \hat{c}_1(m,k) \normtwo{A_k} \hot
\]
yielding (\ref{eq:backerr}) with $c_1(m,k)= 2\sqrt{2}mk + 2\sqrt{k} \geq \sqrt{k} \hat{c}_1 (m,k)$.

To prove (\ref{eq:normaleqbnd}), note that
\[
E_k = \normal{R_k} - \normal{A_k} = \bordermatrix{ & k-1 & 1 \cr k-1 & E_{k-1} & \mbf{w}_k \cr 1 & \mbf{w}_k^T & e_{kk} }
\]
where using Lemma \ref{lem:rkk}, we have
\begin{eqnarray*}
\mbf{w}_k &=& R_{k-1}^T \mbf{s}_k -A_{k-1}^T \mbf{a}_k,\\
e_{kk}&=& \mbf{s}_k^T \mbf{s}_k +r_{kk}^2 -\mbf{a}_k^T \mbf{a}_k \\
&=& \delta_k \mbf{a}_k^T \mbf{a}_k - \Delta_k \mbf{s}_k^T
\mbf{s}_k.
\end{eqnarray*}

Using the bounds on $\delta_k$ and $\Delta_k$ in (\ref{eq:point9}), we have
\begin{eqnarray*}
|e_{kk}| &\leq&  |\delta_k | \normtwo{\mbf{a}_k}^2 + |\Delta_k | \normtwo{\mbf{s}_k}^2 \\
&\leq& (|\delta_k| + |\Delta_k |) \normtwo{\mbf{a}_k}^2 \hot \\
&\leq&  2(m+8) \normtwo{\mbf{a}_k}^2 \macheps \hot \\
&\leq& 2(m+8) \normtwo{A_k}^2 \macheps \hot.
\end{eqnarray*}

Since
\[
\mbf{s}_k + \delta \mbf{s}_k = Q_{k-1}^T \mbf{a}_k, \quad A_{k-1} + \Delta A_{k-1} = Q_{k-1} R_{k-1}
\]
we have
\begin{eqnarray}
\mbf{w}_k &=&  R_{k-1}^T \mbf{s}_k -A_{k-1}^T \mbf{a}_k \nonumber
\\
&=& R_{k-1}^T Q_{k-1}^T \mbf{a}_k -R_{k-1}^T \delta \mbf{s}_k -A_{k-1}^T \mbf{a}_k
\nonumber
\\
&=& \Delta A_{k-1}^T \mbf{a}_k -R_{k-1}^T \delta \mbf{s}_k.
\label{eq:wkexp}
\end{eqnarray}
So  that $\normtwo{\mbf{w}_k}$ has the bound
\begin{eqnarray}
\normtwo{\mbf{w}_k} &\leq& \normtwo{\Delta A_{k-1}}
\normtwo{\mbf{a}_k} + \normtwo{R_{k-1}} \normtwo{\delta
\mbf{s}_k}\hot
\nonumber\\
&\leq&( c_1(m,k-1)\normtwo{A_{k-1}} \normtwo{\mbf{a}_k}+
\sqrt{2(k-1)}m \normtwo{A_{k-1}}\normtwo{\mbf{a}_k})\macheps
\nonumber
\\
&\leq& [2\sqrt{2}m (k-1)+2 \sqrt{k-1} +\sqrt{2(k-1)}m
]\normtwo{A_{k-1}} \normtwo{\mbf{a}_k}\macheps \hot \nonumber
\\
&\leq& 7 m (k-1) \normtwo{A_k}^2 \macheps \hot \label{eq:wkbnd}
\end{eqnarray}

We have that
\begin{eqnarray}
\normtwo{E_k} &\leq& \normtwo{\bear{cc} E_{k-1} & 0 \\ 0 & e_{kk} \enar} + \normtwo{\bear{cc} 0 & \mbf{w}_k \\ \mbf{w}_k^T & 0 \enar}\nonumber \\
 &\leq & \max \{\normtwo{E_{k-1}}, |e_{kk}|\} + \normtwo{\mbf{w}_k}
\nonumber
\\
&\leq& [\max \{ c_2(m,k-1), 2(m+8) \} + 7 m (k-1)]\normtwo{A_k}^2
\macheps\hot \nonumber
\\
&<& [ c_2(m,k-1) + 2(m+8)+ 7 m (k-1) ]\normtwo{A_k}^2 \macheps
\hot \nonumber
\\
&\leq & c_2(m,k) \normtwo{A_k}^2 \macheps \hot
\label{eq:normEk2bnd}
\end{eqnarray}
where
\begin{eqnarray*}
c_2(m,k)&=& \sum_{j=1}^k [2 (m+8) + 7 m (j-1)] \\
&=& 3.5 m (k-1) k + 2mk+16k.
\end{eqnarray*}
Thus we have the expression for $c_2(m,k)$ given in equation (\ref{eq:c2mk}).

To prove (\ref{eq:normRkbnd})--(\ref{eq:normQkbnd}), we simply apply (\ref{eq:backerr})--(\ref{eq:normaleqbnd}).
Equation (\ref{eq:normRkbnd}) results from noting that
\begin{eqnarray*}
\normtwo{R_k}^2 &=& \normtwo{\normal{R_k}} =
\normtwo{\normal{A_k}+E_k}\\
 &\leq&
\normtwo{\normal{A_k}}+\normtwo{E_k} \leq
(1+c_2(m,k)\macheps)\normtwo{A_k}^2 \hot.
\end{eqnarray*}
 Thus,
\[
\normtwo{R_k} \leq (1+c_3(m,k) \macheps) \normtwo{A_k} \hot
\]
where
\[
1+c_3(m,k)\macheps \hot  = \sqrt{1+ c_2(m,k)},
\]
that is, $c_3(m,k)=0.5 c_2(m,k)$. Reversing the roles of $R_k$ and $A_k$ yields
\[
\normtwo{A_k} \leq (1+c_3(m,k) \macheps) \normtwo{R_k} \hot,
\]
thus we have (\ref{eq:normRkbnd}).

To get (\ref{eq:Qkbnd}), we note that
\[
Q_k = (A_k +\Delta A_k)R_k^{-1}
\]
so that
\begin{eqnarray*}
I-\normal{Q_k} &=& R_k^{-T}(\normal{R_k}-\normal{(A_k + \Delta A_k)} ) R_k^{-1}\\
&=& R_k^{-T}(E_k - A_k^T \Delta A_k - (\Delta A_k)^T A_k -
\normal{(\Delta A_k)} ) R_k^{-1}.
\end{eqnarray*}
Thus
\begin{eqnarray*}
\normtwo{I-\normal{Q_k}} &\leq& \normtwo{R_k^{-1}}^2 (\normtwo{E_k} + 2\normtwo{\Delta A_k} \normtwo{A_k} + \normtwo{\Delta A_k}^2)\\
&\leq& \normtwo{R_k^{-1}}^2(c_2(m,k)\normtwo{A_k}^2 + 2c_1(m,k)\normtwo{A_k}^2 + \macheps c_1^2(m,k)\normtwo{A_k}^2)\macheps \hot \\
&\leq& \normtwo{R_k}^2 \normtwo{R_k^{-1}}^2 (c_2(m,k)+2c_1(m,k)) \macheps \hot \\
&=& c_4(m,k) \normtwo{R_k}^2 \normtwo{R_k^{-1}}^2 \macheps \hot
\end{eqnarray*}
where $c_4(m,k) = c_2(m,k)+2c_1(m,k)$.

Finally, to get (\ref{eq:normQkbnd}), we have that
\begin{eqnarray*}
\normtwo{Q_k}^2 &=&\normtwo{\normal{Q_k}} = \normtwo{I-\normal{Q_k}-I}\\
 &\leq& \normtwo{I} + \normtwo{I-\normal{Q_k}} \\
&\leq& 1 + \normtwo{I-\normal{Q_k}} \\
&\leq& 1+ c_4(m,k) \normtwo{R_k}^2\normtwo{R_k^{-1}}^2\macheps \hot
\leq 2 \hot.
\end{eqnarray*}
Taking square roots yields (\ref{eq:normQkbnd}).
\end{proof}

\end{document}